\documentclass[12pt]{article}
\usepackage{amsmath,amsthm, amssymb, floatflt}
\large

\begin{document}

\centerline {\LARGE\bf On S.L.~Tabachnikov's conjecture}
  \medskip

  \centerline {A.I.~Nazarov, F.V.~Petrov}
  \bigskip

  {\small

  S. L. Tabachnikov's conjecture is proved: for any closed curve $\Gamma$
lying inside convex closed planar curve $\Gamma_1$ the mean absolute curvature 
$T(\Gamma)$ exceeds $T(\Gamma_1)$ if $\Gamma\ne k\Gamma_1$.
An inequality $T(\Gamma)\ge T(\Gamma_1)$ is proved for curves 
in a hemisphere.}
  \bigskip

  \section*{1. Problem setting. Main ideas}

Let $\Gamma(s), s\in[0,L(\Gamma)]$ be a naturally parametrized closed curve 
on a plane. We say that $\Gamma(s)$ belongs to the class $BV^1$ if the 
velocity $\Gamma'(s)$ exists and is continious on $[0,L(\Gamma)]$ 
with the exception of countable set; at the 
points of this set $\Gamma'$ has left and right 
limits and the variation of $\Gamma'$ is bounded\footnote
{The variation of a function $f$ mapping into unit circle is defined as 
supremum of sums $\sum_{i=1}^n \rho(f(t_i),f(t_{i-1}))+\rho(f(t_n),f(t_0))$
taken by all subdivisions $t_0<t_1<\dots<t_n$ of a segment in which $f$ is 
defined provided $f$ is defined in the nodes $t_i$; $\rho$ is intrinsic 
metrics of the circle.}. Full variation of $\Gamma'$ is called 
{\it full rotation} of the curve $\Gamma$ and it is denoted by
$V(\Gamma)$.

Note the following properties of the full rotation:

$1^{\circ}$. For $C^2$-smooth curves the full rotation is equal to the 
integral of curvature modulus with respect to natural parameter.

$2^{\circ}$. Full rotation of a closed polygonal line equals the sum of 
external angles in all its vertices.

$3^{\circ}$. Full rotation of a closed convex curve exists and equals $2\pi$.

Define the {\it mean absolute curvature} of a curve $\Gamma\in BV^1$ as its 
full rotation divided by the length: $T(\Gamma)=V(\Gamma)/L(\Gamma)$.

S.L.~Tabachnikov [1] has formulated the following conjecture which he called
{\it DNA inequality}:\medskip

{ \bf Theorem ${\cal P}$.} 1. Mean absolute curvature $T(\Gamma)$ of a closed 
curve $\Gamma\in BV^1$ ("DNA") lying inside convex closed curve $\Gamma_1$
("cell") is not less than $T(\Gamma_1)$.

{2.} If $T(\Gamma)=T(\Gamma_1)$, then curve $\Gamma$ is a multiple circuit of 
$\Gamma_1$.\medskip

A survey of results concerning this conjecture and generalisations is made in 
[1]. The first part of Theorem ${\cal P}$ is proved in [2].

We prove DNA inequality in full generality. The proof of the first part 
partially follows the strategy of [2], but it
is more clear and it is used for
proving the second part. In order to make a paper self-completed,
we give (significantly simplified) proofs of all lemmas from [2] being used.
\medskip

Without loss of generality, we may assume that $\Gamma_1$ is a convex hull 
boundary of $\Gamma$.

A curve $\widetilde\Gamma$ is said {\it better} than a curve $\Gamma$ if 
$T(\Gamma)\ge T(\widetilde\Gamma)$, and it is said {\it strictly better} if 
$T(\Gamma)>T(\widetilde\Gamma)$. We call an {\it improvement} (resp. {strict 
improvement}) of a curve the replacement of a curve onto a better (resp. 
strictly better) one provided the convex hull of "new"\ curve is not larger 
than the convex hull of "old"\ one. Note that if after some improvements of a 
curve $\Gamma$ we get a multiple circuit of $\Gamma_1$ then Claim 1 of 
Theorem ${\cal P}$ is proved for a curve $\Gamma$. If, moreover, at least one 
improvement is strict then the strict inequality $T(\Gamma)>T(\Gamma_1)$ is 
established.

At first Claim 1 is reduced to the case of polygonal lines. 
After that vertices of a polygonal line are moved to a boundary (here and 
further: a boundary of a convex hull). After that each change of rotation 
admits an imrovement of a curve. A finite number of such improvements lets 
us to get a curve, which rotates in only one direction (say, only clockwise), 
for which Claim 1 is almost obvious. Then we prove that every curve 
different from the multiple circuit of a boundary may be strictly improved. 
This proves Claim 2.

\section*{2. Reduction to a polygonal line}

Consider a point of jump of the function $\Gamma'$. Note that the
variation 
of $\Gamma'$ does not change if we redefine $\Gamma'$ at this point as 
arbitrary vector lying on the unit circle between\footnote { i.e. on 
smaller of two arcs.} left and right limits of $\Gamma'$.
Further, speaking about the values set of 
a velocity on some subinterval, we take into account that we add
to it these sets of admissible values in jump points mentoned above. 
So, the values set of velocity on any interval is an arc of unit circle.

We need the following \medskip

{\bf Lemma 1.} Consider two points $A$ and $B$ on a curve $\Gamma$.
A full rotation of curve $\Gamma$ between $A$ and $B$ (such
part of $\Gamma$ will be denoted by $\Gamma_{AB}$) is not less than 
$\rho(\Gamma'(A),e)+\rho(\Gamma'(B),e)$,
where $e$ is a unit vector directed as $\overline{AB}$.
\medskip

{\bf Proof.} If a vector $e$ lies in the value set of $\Gamma'$ at
the part from $A$ to $B$ then the claim is clear
(it suffices to consider a subdivision of
$\Gamma_{AB}$ with nodes $A,\,B$ and a point $C$ with $\Gamma'(C)=e$).
In the opposite case, the values set of $\Gamma'$ from $A$ to $B$ is an 
arc not less than semicircle (otherwise one could find a half-plane
which contains this values set and does not contain $\overline{AB}$). 
Consider a subdivision of $\Gamma_{AB}$, with nodes at points with 
extremal velocities. We get that the rotation of $\Gamma$ from $A$ to 
$B$ is not less than the larger arc between $\Gamma'(A)$ and $\Gamma'(B)$.
This proves Lemma in this case.\hfill$\square$\medskip

{\bf Lemma 2.} Assume that a curve $\Gamma$ does not satisfy Claim 1 
of Theorem ${\cal P}$. Then there exists a polygonal line, which also does 
not satisfy this claim.\medskip

{\bf Proof.} By our assumprions we have $T(\Gamma)<T(\Gamma_1)$.
We inscribe a polygonal line $\Delta$ into $\Gamma$ so that the length 
of $\Delta$ is quite closed to the length of $\Gamma$ (namely, 
$L(\Delta)> L(\Gamma)\cdot \frac {T(\Gamma)}{ T(\Gamma_1)}$).
Clearly, the convex hull of $\Delta$ (by $\Delta_1$ we denote its boundary) 
lies inside $\Gamma_1$.

To prove that $V(\Delta)\le V(\Gamma)$ it suffices to sum up 
the inequalities of Lemma 1 with respect to all edges of a polygonal 
line, and then to use the triangle inequality.

Hence
$$T(\Delta)=\frac {V(\Delta)}{ L(\Delta)}
<\frac {V(\Gamma)}{L(\Gamma)}\cdot \frac
{T(\Gamma_1)}{ T(\Gamma)}=T(\Gamma_1)\le T(\Delta_1),$$

\noindent and Lemma follows. \hfill$\square$\medskip

Let $A_1A_2\dots A_nA_1$ be a closed
polygonal line. We denote by $L$ its length, by $P$
the perimeter of its convex hull and by
$V:=\sum_{i=1}^n(\pi-\angle A_{i-1}A_iA_{i+1})$ its full rotation
(enumeration of indices is cyclic modulo $n$).
We assume that no vertex $A_i$ of a polygonal line
lies on a segment $[A_{i-1}A_{i+1}]$. Such vertices may appear
in a process of improvement, in this case they
will be removed immediately.

In terms of above notations, Claim 1 of Theorem ${\cal P}$
for polygonal lines may be reformulated as follows:\medskip

{\bf Lemma 3.} $\displaystyle{\frac {L}{ V}\le \frac {P}{ 2\pi}}$.\medskip

Note that Lemmas 2 and 3
imply Claim 1 in the general case.

\section*{3. Quadrilaterals}

Here we prove two lemmas, which form a claim of Lemma 3
for quadrilaterals. These lemmas will be used further for improvement
of arbitrary polygonal line.\medskip

{\bf Lemma 4.} For any triangle $ABC$ the inequality 
$$\frac {AB+BC}{ 2\pi-\beta}< \frac {AB+BC+AC}{ 2 \pi},\eqno(1)$$

\noindent  holds with $\beta=\angle ABC$.\medskip

{\bf Proof.} By the sine theorem we have
$$\frac {AB+BC}{ AC}=
  \frac {\sin \angle A+\sin \angle C}{ \sin \beta}=
  \frac {2\sin(\frac {\angle A+\angle C}{ 2})
  \cos(\frac {\angle A-\angle C}{ 2})}{ 2\sin \frac {\beta}{ 2}
  \cos\frac{\beta}{2}}
 =\frac {\cos(\frac {\angle A-\angle C}{ 2})}{ \sin\frac {\beta}{ 2}}\le
\frac {1}{\sin\frac {\beta}{ 2}}.$$

Due to concavity of the sine on $[0,\pi/2]$ we have $\sin \frac {\beta }{ 2}>
  \frac {\beta}{ \pi}$. Hence
$$\frac {AB+BC}{ AC}< \frac {\pi}{ \beta}< \frac {2\pi-\beta}{ \beta},$$

\noindent that is equivalent to (1). \hfill$\square$\medskip

  {\bf Lemma 5.} Let $ABCD$ be a convex quadrilateral, $O$ be a point
of diagonals intersection. Put $\varphi=\angle AOB$. Then
$$\frac {AB+BD+DC+CA}{ 2(\pi+\varphi)}<
\frac {AB+BC+CD+DA}{ 2\pi}.\eqno(2)$$

{\bf Proof.} By Lemma 4 we have $AB> \frac {\varphi}{\pi}(AO+OB)$,
$CD> \frac {\varphi}{\pi}(CO+OD)$. Adding these inequalities we get
$AB+CD> \frac {\varphi}{\pi}(AC+BD)$.

Analogously, $BC+AD>(1-\frac {\varphi}{\pi})(AC+BD)$. Hence
$$\frac {\varphi}{\pi}\cdot \frac {AB+CD}{AC+BD}+
(1+\frac {\varphi}{\pi})\cdot
\frac {BC+DA}{AC+BD}>\frac {\varphi}{\pi}\cdot\frac {\varphi}{\pi}+
(1+\frac {\varphi}{\pi})\cdot(1-\frac {\varphi}{\pi})=1.
$$
Therefore,
$$
\frac {AB+CD} {AC+BD}+1< (1+\frac {\varphi}{\pi})
(\frac {AB+CD}{AC+BD}+\frac {BC+DA}{AC+BD}).
$$
This inequality is equivalent to (2). \hfill$\square$\medskip

Statements of Lemmas 4 and 5 are nothing more than
partial cases of Lemma 3 for concave and selfintersecting
quadrilaterals respectively.
  \medskip

  {\it Remark 1.}
Inequalities of Lemmas 4 and 5 (with $\le$ sign) are true
also for degenerate triangle $ABC$ and quadrilateral $ABCD$ respectively.

\section*{4. Vertices moving to the boundary}

Here we reduce the proof of Lemma 3 to the case,
where all vertices of a polygonal line belong to the boundary
of its convex hull.
\smallskip

Assume that a vertice $A_i$ is situated strictly inside convex hull.
Let us consider three cases.

Case {\bf a)}. The line $A_iA_{i+1}$ does not separate points $A_{i-1}$ and 
$A_{i+2}$. In this case we can strictly improve a polygonal line,
increasing its length without change a rotation:
just move a vertice $A_i$ beyond the segment $A_{i-1}A_i$ while it touches
either the boundary or the ray $A_{i+2}A_{i+1}$.
So we get a better polygonal line
with less number of vertices situated
strictly inside convex hull (it may well be that
a total number of vertices also decreases).
Similar operation is possible if $A_iA_{i-1}$ does not separate $A_{i+1}$ and 
$A_{i-2}$.

Case {\bf b)}. Assume that the line $A_iA_{i+1}$ 
separates $A_{i-1}$ and 
$A_{i+2}$ while $A_iA_{i-1}$ separates $A_{i+1}$ and $A_{i-2}$.

Let $A_{i-2}$ and $A_{i+2}$ be situated in angles supplementary
(by side containing $A_i$) to angles $\angle A_{i+1}A_{i-1}A_i$ and
$\angle A_{i-1}A_{i+1}A_i$, respectively. Then we replace edges $A_{i-1}A_i$ 
and $A_iA_{i+1}$ of a polygonal line $A_1A_2\dots A_n$ to one edge 
$A_{i-1}A_{i+1}$. Assume that the old polygonal line does not satisfy
the inequality of Lemma 3 while new one does, i.e.
$$\frac {L}{ V}> \frac {P}{ 2\pi}\ge
\frac {L-(A_{i-1}A_i+A_iA_{i+1}-A_{i-1}A_{i+1})}{
V-2(\pi-\beta)},\eqno(3)$$
with $\beta=\angle A_{i-1}A_iA_{i+1}$. Then
  $$P\cdot V-2(\pi-\beta)\cdot P+
  2\pi(A_{i-1}A_i+A_iA_{i+1}-A_{i-1}A_{i+1})\ge
  2\pi L> P\cdot V,$$
  hence
$$
 2\pi(A_{i-1}A_i+A_iA_{i+1}-A_{i-1}A_{i+1})
  > 2(\pi-\beta) P
  \ge2(\pi-\beta)(A_{i-1}A_i+A_iA_{i+1}+A_{i-1}A_{i+1}),
$$
and we obtain a contradiction with Lemma 4. Thus, the new polygonal line
also has to be a counterexample to Lemma 3, while it has less inner
vertices.

It remains to consider case {\bf c)}, where, for example, the line
$A_{i-1}A_{i+1}$ separates $A_{i+2}$ and $A_i$ (in this case a vertice 
$A_{i+1}$ also lies strictly inside the convex hull).
Without loss of generality we assume that the angle 
$A_{i-1}A_{i+1}A_i$ is the least for all 
indices $i$ satisfying this condition. Let us replace $i$ to $i+1$, 
and consider analogous cases. The 
vertex $A_{i-1}$ lies in an angle supplementary 
to $\angle A_iA_{i+2}A_{i+1}$ with respect to $A_iA_{i+1}$. If the vertex 
$A_{i+3}$ does not lie in angle vertical to $A_{i+1}A_{i+2}A_i$, the 
polygonal line may be improved as it is shown (with the change 
$i\rightarrow i+1$). In the opposite case we get a contradiction with the 
choise of $i$: the angle $\angle A_iA_{i+2}A_{i+1}$ is less than the angle 
$\angle A_iA_{i+1}A_{i-1}$ (since 
$\angle A_iA_{i+2}A_{i+1}+\angle A_iA_{i+1}A_{i+2}<
\angle A_iA_{i+1}A_{i_1}+\angle A_iA_{i+1}A_{i+2}$).

So, by the finite number of steps we reduce general case to the case,
where all the vertices $A_i$ of a polygonal line $A_1A_2\dots A_n$ lie on the 
boundary of its convex hull.

\section*{5. Decreasing the number of direction changes}

Fix an orientation of the plane. We say that polygonal line $A_1A_2\dots A_n$ 
turns to the right in a vertex $A_i$ if the base $\overline{A_{i-1}A_i}$, 
$\overline{A_iA_{i+1}}$ is negatively oriented. Otherwise
(in particular, if these vectors are collinear) we say that it turns to the 
left.

If the polygonal line turns to the right (or to the left) two times in succession we 
may replace an edge between these turns to the part of boundary, passed
in the same direction. This operation improves a polygonal line. We call it 
{\it the stretching} of a polygonal line.\medskip

{\bf Lemma 6.} Assume that some consequent edges of our polygonal line form a 
full circuit of a boundary, and the first and the last edges coincide (i.e. 
the boundary is a convex polygon $C_1C_2\dots C_m$ while a polygonal line has 
a part $XC_1C_2\dots C_mC_1C_2Y$). Then the claim of Lemma 3 for this 
polygonal line is equiavelent to the claim of Lemma 3 for a polygonal line 
with this part removed (i.e. for a polygonal line, in which this part
is replaced by just $XC_1C_2Y$). \medskip

{\bf Proof.} Note that perimeter of a polygonal line after circuit removing 
equals $L-P$, and its full rotation equals $V-2\pi$. The statement of Lemma 3 
for the new polygonal line claims $\frac{L-P}{V-2\pi}\le \frac{P}{2\pi}$, 
which is equivalent to the statement of Lemma 3 for the
initial polygonal line.
\hfill$\square$\medskip

Let us repeat the operation of Lemma 6 while it is possible. This process 
must stop while the number of edges decreases. Note that the
number of changes of turns 
directions does not change.

Now the polygonal line is partitioned to parts, in which all the turns have 
the same direction; and in each part all the edges (except, possibly, the 
first and the last) go along the boundary and, due to Lemma 6, are distinct.

Let us develop the following operation. Choose a part $A_iA_{i+1}\dots A_k$,
in which all turns are, say, the left ones (namely, turns in vertices 
$A_{i+1}$, $A_{i+2},\dots$, $A_{k-1}$ are left, and turns in $A_i$ and $A_k$ 
are right). 
Replace the path $A_iA_{i+1}\dots A_k$ to the part of boundary $A_i..A_k$ 
bypassing the boundary in the opposite (in our case --- negative) direction. 
The number of direction changes decreases after this operation.

Assume that the initial polygonal line does not satisfy the inequality of Lemma 
3. Our goal is to prove that in this case a new polygonal line does not 
satisfy it as well. Six cases are possible, they are determined by the
order of 
vertices $A_i$, $A_{i+1}$, $A_{k-1}$, $A_k$ while bypassing a boundary in 
the positive direction:

  $1^{\circ}$. $A_iA_{i+1}A_{k-1}A_k$;

$2^{\circ}$. $A_iA_kA_{i+1}A_{k-1}$;

$3^{\circ}$. $A_iA_{i+1}A_kA_{k-1}$;

$4^{\circ}$. $A_iA_{k-1}A_kA_{i+1}$;

$5^{\circ}$. $A_iA_{k-1}A_{i+1}A_k$;

$6^{\circ}$. $A_iA_{k-1}A_kA_{i+1}$.

Cases 3 and 4 are equivalent up to renaming and symmetry. Denote the 
length of a polygonal line $A_iA_{i+1}\dots A_k$ by $s$, and denote
the length of its replacement by $s'$.
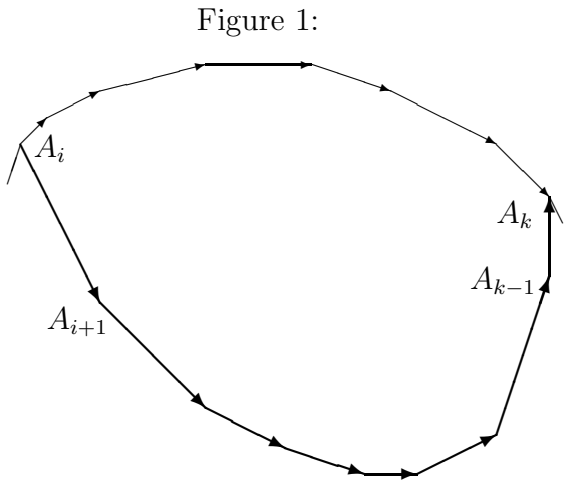
\medskip

\begin{floatingfigure}{240pt}

\caption{}

\begin{picture}(220,200)(0,40)

\put(15,195){$A_i$}
\put(20,130){$A_{i+1}$}
\put(180,145){$A_{k-1}$}
\put(190,170){$A_k$}

\put(10,200){\line(-1,-3){5}}
\put(10,200){\vector(1,1){10}}
\put(20,210){\vector(2,1){20}}
\put(40,220){\vector(4,1){40}}
\put(80,230){\vector(1,0){40}}
\put(120,230){\vector(3,-1){30}}
\put(150,220){\vector(2,-1){40}}
\put(190,200){\vector(1,-1){20}}

\put(210,180){\line(1,-2){5}}

\thicklines

\put(10,200){\vector(1,-2){30}}
\put(40,140){\vector(1,-1){40}}
\put(80,100){\vector(2,-1){30}}
\put(110,85){\vector(3,-1){30}}
\put(140,75){\vector(1,0){20}}
\put(160,75){\vector(2,1){30}}
\put(190,90){\vector(1,3){20}}
\put(210,150){\vector(0,1){30}}

\end{picture}

\end{floatingfigure}

$1^{\circ}$  (see Figure 1). Denote $\angle A_{i+1}A_iA_k=\alpha$, 
$\angle A_{k-1}A_kA_i=\beta$. After the change of a polygonal line, its full 
rotation decreases by $2(\alpha+\beta)$, and the length decreases by $s-s'$.

If the new polygonal line satisfies the inequality of Lemma 3, we have
$$\frac {L}{ V}>\frac {P}{ 2\pi}\ge 
\frac {L+s'-s}{V-2(\alpha+\beta)},\eqno(4)$$
or
$$P\cdot V -2P(\alpha+\beta)+2\pi (s-s')\ge 2\pi L>P\cdot V,$$
hence $2\pi (s-s')>2(\alpha+\beta)P\ge 2(\alpha+\beta)(s+s')$, and
$2(\pi-\alpha-\beta)s>2(\pi+\alpha+\beta)s'$.

The last inequality may hold only if $\alpha+\beta<\pi$. In this case rays
$A_iA_{i+1}$ and $A_kA_{k-1}$ meet in point $C$, and $CA_i+CA_k\ge s$, 
$A_iA_k\le s'$. So, 
$$(\pi-\alpha-\beta)(CA_i+CA_k)>(\pi+\alpha+\beta) A_iA_k,$$
that contradicts Lemma 4.\smallskip

$2^{\circ}$ (see Figure 2). Denote by $O$ the point of intersection of segments 
$A_iA_{i+1}$ and $A_kA_{k-1}$, $\angle A_iOA_k=\varphi$. After the change of 
a polygonal line, its full rotation decreases by  $2\varphi$, and the length 
decreases by $s-s'$.

\begin{floatingfigure}{240pt}

\caption{}

\begin{picture}(220,220)(0,50)

\put(130,235){$A_i$}
\put(20,130){$A_{i+1}$}
\put(180,145){$A_{k-1}$}
\put(10,215){$A_k$}

\put(30,210){\line(2,1){20}}
\put(130,230){\vector(2,-1){50}}
\put(180,205){\vector(1,-1){20}}
\put(200,185){\vector(1,-2){10}}
\put(210,165){\vector(0,-1){15}}
\put(40,140){\vector(-1,2){20}}
\put(20,180){\vector(0,1){20}}
\put(20,200){\vector(1,1){10}}

\put(130,230){\line(-4,1){20}}

\thicklines

\put(130,230){\vector(-1,-1){90}}
\put(40,140){\vector(1,-1){40}}
\put(80,100){\vector(2,-1){30}}
\put(110,85){\vector(3,-1){30}}
\put(140,75){\vector(1,0){20}}
\put(160,75){\vector(2,1){30}}
\put(190,90){\vector(1,3){20}}
\put(210,150){\vector(-3,1){180}}

\end{picture}

\end{floatingfigure}

If new polygonal line satisfies Lemma 3 then
$$\frac {L}{V}>\frac {P}{2\pi}\ge \frac {L+s'-s}{V-2\varphi}.$$
From here, analogously to p.$1^{\circ}$, we obtain
$2\pi (s-s')>2\varphi P$, whence
$$\frac {A_iA_{i+1}+A_kA_{k-1}-(A_iA_{k-1}+A_kA_{i+1})}
{A_iA_k+A_kA_{i+1}+A_{i+1}A_{k-1}+A_iA_{k-1}}>\frac {\varphi}{\pi}.$$
This contradicts Lemma 5 (for quadrilateral $A_iA_{k-1}A_{i+1}A_k$).

Other cases are analogous to these two, Lemma 4 is used in cases 3 and 6 and 
Lemma 5 is used in case 5.

So, after a finite number of steps we get a polygonal line which turns only 
to the right. Using stretching we get a multiple circuit of a boundary from 
this polygonal line, hence the statement of Lemma 3 holds for this
polygonal line. So, initial assumption was wrong, and Lemma
3 is proved. Claim 1 of Theorem $\cal P$ is proved as well.

\section*{6. Proof of Claim 2}

Assume that a curve $\Gamma$ is not a (multiple) boundary circuit, but
$T(\Gamma)=T(\Gamma_1)$. We select a finite number of points on $\Gamma$ so 
that sum of velocity jumps in other points is quite small (say, less then 
$\pi/180$). A union of this finite set and a set $\Gamma\cap\Gamma_1$ is 
closed. Preimage (recall that a curve is naturally parametrised) of its 
complement is a union of countable number of intervals. Consider one of these 
intervals, let it corresponds to the part of $\Gamma$ between points $A$ and 
$B$.

A part $\Gamma_{CD}$ is said to be {\it small} if the values set of velocity 
on this part is an arc of length at most $\pi/4$, and a circle with diameter
$CD$ lies strictly inside $\Gamma_1$. It is easy to see that any inner point 
of $\Gamma_{AB}$ belongs to some small subpart.

Consider a small part $\Gamma_{CD}$. Redefine, if necessary, velocities
$\Gamma'(C)$ and $\Gamma'(D)$ as their right limits. Define a parallelogram 
$CPQD$ with vectors $\overline{CP}$ and $\overline{CQ}$ directed as extremal
directions of the velocity of curve $\Gamma$ on a part $\Gamma_{CD}$. This 
parallelogram lies strictly inside $\Gamma_1$ (since extremal directions are 
quite close to the direction of vector $\overline{CD}$).

There are points $X$ and $Y$ on $\Gamma_{CD}$ with tangents parallel to $CP$ 
and $CQ$, respectively. Without loss of generality we can say
that the order of points is 
$C-X-Y-D$. Replace a part $\Gamma_{CD}$ to a polygonal line $CPD$.

Note that full rotation of $\Gamma_{CD}$ is not less than
$$v:=\rho(\Gamma'(C),\Gamma'(X))+\rho(\Gamma'(X),\Gamma'(Y))+
\rho(\Gamma'(Y),\Gamma'(D)),$$
while the full rotation of a new part equals $v$. Equality holds only if 
$\Gamma$ is convex from $C$ to $X$, from $X$ to $Y$ and from $Y$ to $D$.

Furthermore, the length of $\Gamma_{CD}$ does not exceed $CP+PD$. To prove 
this consider an arbitrary polygonal line inscribed in $\Gamma_{CD}$. Directions 
of its edges may differ between directions of vectors $\overline{CP}$ and 
$\overline{CQ}$. So, after their rearrangement in the order monotone (in the 
sense of direction) from $\overline{CP}$ till $\overline{CQ}$ we get a
convex polygonal line $C\dots D$ situated inside triangle $CPD$, and hence 
its length does not exceed $CP+PD$.

So, after replacement of a part $\Gamma_{CD}$ onto a polygonal line
$CPD$ the full rotation does not decrease, and length does not increase, i.e. 
curve $\Gamma$ improves. But it may not be strictly improved due to Claim 1 
already proved. Hence after such replacement neither full rotation,
nor length changes. The first is possible only if $\Gamma_{CD}$ may be 
splitted into at most three convex parts ($C-X$, $X-Y$, $Y-D$). The second is 
possible only if each of these parts is a polygonal line with at most two 
edges. So, $\Gamma_{CD}$ is a polygonal line with at most six edges
(the number of edges may be decreased, but we do not need it).\smallskip

Now we fix points $A'$ and  $B'$ on an open arc $\Gamma_{AB}$ and cover 
$\Gamma_{A'B'}$ by finite number of small parts. Curve $\Gamma$ is a 
polygonal line on each small part, hence $\Gamma_{A'B'}$ is a polygonal line 
too.

Now we prove that if $\Gamma_{A'B'}$ has at least four edges then $\Gamma$ 
may be strictly improved. Consider cases of \S4. In the case {\bf a)} we use 
only local structure of a polygonal line, and the same argument works in our 
situation.

In the case {\bf b)}, using Claim 1 for a changed curve, we obtain the
inequality analogous to (3):
$$\frac {L(\Gamma)}{V(\Gamma)}=\frac {L(\Gamma_1)}{2\pi}\ge
\frac {L(\Gamma)-(A_{i-1}A_i+A_iA_{i+1}-A_{i-1}A_{i+1})}
{V(\Gamma)-2(\pi-\beta)},$$

\noindent hence
$$2\pi(A_{i-1}A_i+A_iA_{i+1}-A_{i-1}A_{i+1}) \ge
  2(\pi-\beta)(A_{i-1}A_i+A_iA_{i+1}+A_{i-1}A_{i+1}),$$
that contradicts Lemma 4.

In the case {\bf c)}, if $A_{i+2}$ and $A_i$ are separated by line 
$A_{i-1}A_{i+1}$ (in this case a vertex $A_{i+1}$ also lies strictly
inside convex curve), the polygonal line may be strictly improved by 
replacing an edge $A_iA_{i+1}$ to a parallel longer edge
$A_i'A_{i+1}'\parallel A_iA_{i+1}$, where
$A_i\in [A_{i-1}A_i'[$, $A_{i+1}'\in ]A_{i+1}A_{i+2}]$.

So, $\Gamma_{A'B'}$ is a polygonal line with at most three edges. Since
$A'$ and $B'$ were chosen arbitrary, a curve $\Gamma_{AB}$ is a polygonal 
line with at most three edges too.\smallskip

Now we add the points of "large turn", excluded before, to the considered
intervals. Then the whole inner part of a curve $\Gamma$ is splitted to at 
most countable set of polygonal lines with a finite number of edges.

If one of such parts contains more than one edge, the curve $\Gamma$ may be 
strictly improved as it was done in \S4.\smallskip

So, all points of $\Gamma$ lie either on a boundary, or on segments
joining boundary points.

Assume that the number of segments is infinite. Then there exists a sequence 
of segments with length tending to 0 and endpoints tending to some point 
$C\in\Gamma_1$. Let us fix a small neighborhood of the point $C$ with full 
rotation of $\Gamma_1$ equal to $\varphi_0<\pi$. Consider one of segments
${\overline{AB}}$ ($A,\,B\in \Gamma_1$) lying in this neighborhood.
If vectors $\Gamma'(A-)$ and $\Gamma'(B+)$ are directed to different sides 
with respect to the line $AB$, then the curve $\Gamma$ may be strictly 
improved by the stretching of a segment $AB$ to the boundary. It is 
impossible. So, the variation of $\Gamma'$ on $AB$ is not less than 
$\pi-\varphi_0$, hence the full variation is infinite.\smallskip

So, the curve $\Gamma$ consists of a finite number of boundary pieces
and a finite number of segments between them. If $\Gamma$ contains return 
points on a boundary (since $\Gamma\in BV^1$ there may be only a 
finite number of such points), we consider them as 
"inner segments of zero length".

If two consecutive pieces of the boundary have the same circuit direction,
it admits an improvement of $\Gamma$: just stretch the segment between them.
Further, we may remove all the full circuits of a boundary as in Lemma 6.

Now consider an arc $\Gamma_{AB}$ consisting of the segment $AA_1$, piece of 
boundary $\Gamma_{A_1B_1}$ (which has, say, positive direction), and the 
segment $B_1B$. Analogously to \S5, replace $\Gamma_{AB}$ to a "negative" arc 
of a boundary between $A$ and $B$. Here we have to consider six cases of
\S5 again, cases depend on the order of points $A$, $A_1$, $B$, $B_1$ in a 
positive circuit. For example, in the case $1^{\circ}$ (order $AA_1B_1B$), we 
use Claim 1 for the new curve and get the inequality analogous to (4),
  $$\frac {L(\Gamma)}{V(\Gamma)}=\frac {L(\Gamma_1)}{ 2\pi}\ge \frac
{L(\Gamma)+s'-s}{V(\Gamma)-2(\angle A_1AB+\angle B_1BA)}.$$

From here we deduce, as in \S5, that the rays $AA_1$ and $BB_1$
meet in a point $C$, and
$$(\pi-\angle A_1AB-\angle B_1BA)\cdot(CA+CB)\ge
(\pi+\angle A_1AB+\angle B_1BA)\cdot AB,$$
that contradicts Lemma 4.\smallskip

Analogous contradictions can be obtained in remaining cases. It shows that 
the curve $\Gamma$ may not have inner segments, and hence it is a circuit
of a boundary. Now remember the circuits removed earlier and realize that 
initially $\Gamma$ was a multiple circuit of a boundary. The statement 2 is 
proved.\hfill$\square$\medskip

\section*{7. The surfaces of constant curvature}

Now we prove a statement generalizing DNA inequality to a spherical case.

Let $\Gamma$ be a closed curve lying in some hemisphere (here and further: 
of unit radius). Let the variation of the right rotation $V(\Gamma)$ be 
finite. For the definitions
we refer to [3]. Note that if $\Gamma:\ A_1A_2\dots A_nA_1$ is 
a closed polygonal line then 
$V(\Gamma)=\sum_{i=1}^n (\pi-\angle A_{i-1}A_iA_{i+1})$ (enumeration of 
indices is cyclic).

Define a {\it mean absolute geodesic curvature} $T(\Gamma)$ of a closed curve 
on a sphere as $T(\Gamma)=V(\Gamma)/L(\Gamma)$.
\medskip

{ \bf Theorem ${\cal S}$.}
Let $\Gamma$ be a closed curve in a hemisphere, and let the variation 
of its right rotation be finite. Let $\Gamma_1$ be a boundary of its convex 
hull.
Then 
$T(\Gamma)\ge T(\Gamma_1)=(2\pi-S)/L(\Gamma_1)$, where $S$ is the area of the 
convex hull.\medskip

The plan of a proof of a Theorem ${\cal S}$ is the same as in planar case. 
First of all, we formulate corresponding statement for poligonal lines.\medskip

{ \bf Theorem ${\cal S}'$.} Let $\Gamma$ be a closed poligonal line in a 
hemisphere, and let $\Gamma_1$ be a boundary of its convex hull.
If $\Gamma$ is not  multiple circuit of $\Gamma_1$ then $T(\Gamma)> 
T(\Gamma_1)$.\medskip

Before we pass to the case of quadrilaterals, we prove the following claim 
elaborating (in particular case) the theorem of A.D.~Aleksandrov on angles 
comparing. \medskip

{\bf Lemma $1s$.} Let $ABC$ be a non-degenerate triangle on a sphere. We
denote its sides by $BC=a$, $CA=b$, $AB=c$, and its angles by
$\alpha$, $\beta$, $\gamma$ respectively. 
We denote by $\alpha'$, $\beta'$, $\gamma'$ the angles
of a triangle with sides $a$, $b$, $c$ on a plane. Then
$$\alpha-\alpha'< (\beta-\beta')+(\gamma-\gamma'). \eqno(1_s)$$

{\bf Proof.} We denote $a+b+c=4S$, $S-a/2=X$, $S-b/2=Y$, $S-c/2=Z$.
${\cal E}=\alpha+\beta+\gamma-\pi$ is an area of triangle $ABC$.
The inequality $(1_s)$ is equivalent to the inequality
$\alpha'> \alpha-{\cal E}/2$ or, in other words,
 $$\tan \frac{\alpha'}2 > \tan(\alpha/2-{\cal E}/4).\eqno(2_s)$$
Substituting here the formulas
$$\gathered
\tan\frac {\alpha} 2=\sqrt\frac {\sin 2Y \sin 2Z} {\sin 2X \sin 2S},\quad
\tan\frac {\cal E}4=
\sqrt {\tan S\cdot \tan X\cdot\tan Y\cdot\tan Z },\\
\tan \frac {\alpha'}2= \sqrt {\frac{YZ}{X(X+Y+Z)}}
\endgathered$$
 (the first formula is [3, (28)], the second is [4],
 the third is [5, (20)]),
we convert ($2_s$) to the inequality
 $$
\frac {\sin Z\sin Y} {\sin X\sin S}\cdot
\frac {\cos Y\cos Z-\sin X\sin S}{\cos X\cos S+\sin Y\sin Z}<
 \sqrt{\frac {YZ} {X(X+Y+Z)}}\cdot\sqrt {\frac {\sin(2Z)\sin(2Y)}
 {\sin(2X)\sin(2S)}}.\eqno (3_s)
$$
 Since $S=X+Y+Z$, we have $\cos(S-X)=\cos(Y+Z)$,
 hence the second multiple in the left-hand side of ($3_s$) equals 1.
Let us denote $f(x)=x\cot x$, then ($3_s$) reduces to
 $$f(Y)f(Z)> f(X)f(X+Y+Z). \eqno(4_s)$$
Since $f'(x)=\frac {sin (2x)-2x}{2\sin^2 x}<0$  for $0<x<\frac{\pi}2$,
the function $f$ srictly decreases on 
$[0,\frac{\pi}2]$. Since all the arguments in 
($4_s$) lie in $[0,\frac{\pi}2]$ (we recall that $X+Y+Z=(a+b+c)/4\le \pi/2$),
we may suppose that $X=0$ and prove the inequality
 $$f(Y)f(Z)> f(0)f(Y+Z).\eqno(5_s)$$
We have $(\ln (f))''(x)= \frac {4}{\sin^2(2x)}\left(\cos(2x)-
\frac {\sin^2(2x)}{4x^2}\right)$. We omit an elementary proof of the
inequality
$\cos t< (\frac {\sin t}t)^2$ for $t=2x\in]0,\pi[$. This shows that 
$\ln (f)$ is strictly concave on $[0,\frac{\pi}2]$, and ($5_s$) 
follows.\hfill$\square$\medskip

Now we are ready to prove the analogs of Lemmas 4 and 5 on a sphere.
\medskip

{\bf Lemma $2s$}. Let $ABC$ be a non-degenerate triangle on a sphere.
Then, with the same notations as in Lemma $1s$,
$$
\frac {a+c}{2\pi-\beta}< \frac {a+b+c}{2\pi-{\cal E}}.
$$

{\bf Proof.} The statement follows from a chain of inequalities
$$\frac {a+c}{a+b+c}<\frac {2\pi-\beta'}{2\pi}<
\frac {2\pi-\beta+{\cal E}/2}{2\pi}\le \frac
{2\pi-\beta}{2\pi-{\cal E}}$$
(the first inequality is Lemma 4, the second is Lemma 1s,
the third reduces to the obvious $\beta\le \pi+{\cal E}/2$).
\hfill$\square$\medskip

{\bf Lemma $3s$}. 
 Consider a convex spherical quadrilateral $ABCD$ on a hemisphere
and denote by 
$O$ the point of diagonals intersection. Put $\varphi=\angle AOB$. Denote 
$AB=a$, $BC=b$, $CD=c$, $DA=d$, $BD=m$, $AC=n$, and $\angle AOB=\varphi$.
We denote by ${\cal E}_1$, ${\cal E}_2$, ${\cal E}_3$, ${\cal E}_4$ 
the areas of 
triangles $OAB$, $OBC$, $OCD$, $ODA$, respectively, and put 
${\cal E}={\cal E}_1+{\cal E}_2+{\cal E}_3+{\cal E}_4$.
Then
$$
\frac {a+c+m+n}{2\pi-({\cal E}_1+{\cal E}_3)+2\varphi}<
\frac {a+b+c+d}{2\pi-{\cal E}}.\eqno(6_s)
$$

{\bf Proof.} We denote by $\varphi'$ an angle of a planar triangle with
sides $a$, $AO$, $BO$, opposite to side $a$. Then
$$\frac a{AO+BO}>\frac {\varphi'}{\pi}> \frac {\varphi-{\cal E}_1/2}{\pi}
>\frac {\varphi-({\cal E}_1+{\cal E}_3)/2}{\pi}, \eqno (7_s)$$
(the first inequality follows from the proof of Lemma 4, the second one --
from the Lemma 1s). Analogously
$$\frac c{CO+DO}>\frac {\varphi-({\cal E}_1+{\cal E}_3)/2}{\pi}. \eqno (8_s)$$
Estimates ($7_s$) and ($8_s$) imply that
$$x:=\frac {a+c}{m+n}>\frac {\varphi-({\cal E}_1+{\cal E}_3)/2}{\pi}. 
$$
Analogously,
$$
y:=\frac {b+d}{m+n}>\frac {\pi-\varphi-({\cal E}_2+{\cal E}_4)/2}{\pi}. 
$$

We substitute the lower bounds for $x$ and $y$ to the equality
$$
z:=\frac {a+c+m+n}{a+c+b+d}=\frac {x+1}{x+y}=1+\frac {1-y}{x+y}.
$$
Since $y<1$, we get an upper bound for $z$. 
It gives $(6_s)$.\hfill$\square$\medskip

Now we briefly explain the plan of the proof of Theorem ${\cal S}'$.

Arguments used in cases {\bf a)} and {\bf c)} of Section 4 may be transfered 
with minor changes (natural changes arising due to spherical excess $\cal E$ 
play for us).

In the case {\bf b)} we choose index $i$ so that an angle vertical to 
$\angle A_{i-1}A_{i+1}A_i$ has the least area in hemisphere.
Angle vertical to $\angle A_iA_{i+2}A_{i+1}$ is contained in an angle 
vertical to $\angle A_{i-1}A_{i+1}A_i$ (in hemisphere),
since there are no conjugate points on a hemisphere. This leads to a 
contradiction.

The arguments of Section 5 are changed in the same way as ones of case 
{\bf c)} of Section 4.

{ \bf The deduction of Theorem ${\cal S}$ from 
Theorem ${\cal S}'$.} First of all we note
that the curve $\Gamma$ can be splitted into a 
finite number of parts without
self-intersections. Really, if the part of $\Gamma$ with self-intersections 
has sufficiently small length then its rotation is not less than $\pi/2$. 

Let $0=t_1<t_2<\dots<t_n<t_{n+1}=L(\Gamma)$ be the nodes of this partition
(we can suppose $\Gamma$ is naturally parametrized with the starting point 
in a node). Set $A_i:=\Gamma(t_i)$, $\Gamma_i:=\Gamma_{[t_i,t_{i+1}]}$. By 
Theorem 1 [4], for $i=1,\dots,n$ there exists a sequence of poligonal lines
$g_k^i:\ [t_i,t_{i+1}]\rightarrow S^2$ $(k=1,\,2,\dots)$ s.t. 
$g_k^i(t_i)=A_i$, $g_k^i(t_{i+1})=A_{i+1}$, $g_k^i$ converge to $\Gamma_i$ 
from the right, and $\limsup V(g_k^i)\le V(\Gamma_i)$. Moreover, the 
directions of $g_k^i$ at the points $A_i$ and $A_{i+1}$ converge to the 
directions (right and left, correspondingly) of the curve $\Gamma$ at these 
points.

We denote by $g_k$ the sum of poligonal
lines $g_k^i$ with respect to $i=1,\dots,n$, 
and by $G_k$ the boundary of the convex hull of $g_k$. Then
$\limsup V(g_k)\le V(\Gamma)$. Further, $L(\Gamma)\le \liminf L(g_k)$. 
Finally, since $g_k\rightarrow \Gamma$ uniformly, then 
$G_k\rightarrow \Gamma_1$, whence $L(G_k)\rightarrow L(\Gamma_1)$ and 
$T(G_k)\rightarrow T(\Gamma_1)$ . By Theorem ${\cal S}'$
$$T(\Gamma)\ge \limsup T(g_k)\ge \limsup T(G_k)= T(\Gamma_1).$$
This completes the proof.\hfill$\square$\medskip

Unfortunately, we cannot transfer the statement 2 of Theorem $\cal P$ to the 
sphere.\bigskip

Note that in the Lobachevskii plane the DNA inequality is not true. To show 
this we consider 
in the Lobachevskii plane a triangle $ABC$ and a polygonal 
line $\Gamma=ABCC_1B_1BCA$ with some points $B_1\in AB$, $C_1\in AC$. Then 
$\Gamma_1=ABCA$, and 
$$\frac {V(\Gamma)}{V(\Gamma_1)}=2-\frac {S(A_1B_1C_1)}{V(\Gamma_1)}<2-\frac 
{S(A_1B_1C_1)}{3\pi}.$$ 
Moving the vertex $B$ sufficiently far along the ray $AB_1$ we make the 
quotient $L(\Gamma)/L(G_1)$ arbitrary close to $2$ that gives
$T(\Gamma)<T(\Gamma_1)$.\bigskip

We are grateful to V.A.~Zalgaller who attracted our attention to the 
references [3], [4]. We also thank S.V.~Duzhin for his attention.
\medskip

The paper is supported by grant 
VNSh-2261.2003.1 (the first author) and 
by grants VNSh-2251.2003.1 and RFFR 02-01-00093 (the second author).
 \bigskip

  {\centerline {References} }
\medskip

1. S. Tabachnikov. The tale of a geometric inequality // MASS colloquium
lecture, 2001.

2. J. Lagarias, T. Richardson. Convexity and the average curvature of
the plane curves // Geom. Dedicata, 67 (1997), 1-38.

3. A.D. Aleksandrov. Intrinsic Geometry of Convex 
Surfaces. OGIZ, Moscow --
Leningrad, 1948 (Russian). English transl.: A.D. Alexandrov. 
Selected works: 
Intrinsic Geometry of Convex Surfaces. CRC, 2005.

4. V. A. Zalgaller. On curves with curvature of bounded 
variation on a 
convex surface // Mat. Sbornik, 26 (1950). 205--214 (Russian).

5. http:$\backslash\backslash$ mathworld.wolfram.com$\backslash$
SphericalTrigonometry.html 

6. http:$\backslash\backslash$ mathworld.wolfram.com$\backslash$
SphericalExcess.html

7. http:$\backslash\backslash$ mathworld.wolfram.com$\backslash$
Triangle.html

\end{document}